\documentclass{amsart}
\usepackage{amsfonts,amssymb,amsmath,amsthm}
\usepackage{url}
\usepackage{enumerate}
\newtheorem{theorem}{Theorem}[section]
\newtheorem{proposition}[theorem]{Proposition}

\newtheorem{corollary}[theorem]{Corollary}
\theoremstyle{definition}
\newtheorem{definition}[theorem]{Definition}
\newtheorem{example}[theorem]{Example}

\newtheorem{remark}[theorem]{Remark}

\numberwithin{equation}{section}

\begin{document}

\title{Locally convex inductive limit cones}

\author{M.R. Motallebi}
\address{Department of Mathematics, Faculty of Mathematical Sciences, University of Mohaghegh  Ardabili, P.O. Box 179,  Ardabil, Iran}
\email{motallebi@uma.ac.ir}



\begin{abstract}
We define the finest order on inductive limits of ordered cones
which makes the linear mappings monotone and
 gives rise to
the definition of inductive limit topologies for cones. Using the
polars of neighborhoods, we establish  embeddings between  direct
sums, inductive limits and their duals.  These lead us to
investigate the weak topologies and the topologies of pointwise
convergence in inductive limits.
\end{abstract} \maketitle

\section{Introduction} \label{sect1}
The theory of locally convex cones constitute canonical
generalizations of locally convex topological vector spaces, still
retaining many of their most important properties. In particular,
locally convex cones yield rich duality theory, where the study of
the dual cone offers valuable insight into the given locally
convex cone itself. The locally convex cone topologies can be
generated using dual pairs of cones and a bilinear form which
leads to various notions of weak and polar topologies. There are
the notions of continuous linear mappings between locally convex
cones and their adjoints as linear mappings between the respective
duals (see [1-6, 9]).
 The topology of a locally convex cone is given
through a convex quasi-uniform structure, but it can alternatively
be expressed in terms of its order structure alone. The aim of this
paper is the study of  inductive limit topology  and its duality
properties.
 We introduce inductive limit order of a family of
ordered cones  by the linear mappings, i.e., the finest order such
that the linear mappings  are monotone. This uses the product of
abstract neighborhoods
 to define the inductive limit topologies  for cones.
We discuss the linear mappings  and the polars of neighborhoods
for inductive limits; in particular,
 we show that an inductive limit topology  will carry  the strict separation property,
 whenever its components have this property.
Using the adjoints of linear mappings,
 we embed  direct sums on inductive limits and inductive limit dual cones in direct sum dual
 cones which lead  to the study of  weak topologies and the
 topologies of pointwise convergence for inductive limits.
   \par A \emph{cone} is a set $\mathcal{P}$ endowed with an addition
$(a,b)\longmapsto a+b$ and scalar multiplication $(\alpha,
a)\longmapsto \alpha a$ for real numbers $\alpha\geq 0.$
 The addition is supposed to be associative and commutative, there
is a neutral element $0\in \mathcal{P}.$ For the scalar
multiplication the usual associative and distributive properties
hold, that is, $\alpha (\beta
a)=(\alpha\beta)a,\,\,(\alpha+\beta)a=\alpha a+\beta
a,\,\,\alpha(a+b)=\alpha a+\alpha b,\,\,1 a=a,\,\,0 a=0$ for all
$a,b\in\mathcal{P}$ and $\alpha,\beta\geq 0.$
 An \textit{ordered cone}
$\mathcal{P}$ carries  a reflexive transitive relation $\leq$ such
that $a\leq b$ implies $a+c\leq b+c$ and $\alpha a \leq\alpha b$ for
all $a,b,c\in \mathcal{P}$ and $\alpha\geq 0.$ Equality is such an
order, hence cones without an explicit order structure are also
included. Note that anti-symmetry is note required for the relation
$\leq.$ For example, the extended scalar field
$\overline{\mathbb{R}}=\mathbb{R}\cup\{+\infty\}$ of real numbers is
an ordered cone. We consider the usual order and algebraic
operations in $\overline{\mathbb{R}};$ in particular,
$\alpha+\infty=+\infty$ for all $\alpha\in\overline{\mathbb{R}},$
$\alpha \cdot(+\infty)=+\infty$ for all $\alpha>0$ and $0\cdot
(+\infty)=0.$
\par A \emph{full locally convex cone} $(\mathcal{P},\mathcal{V})$ is an ordered cone $\mathcal{P}$
 that contains
 an \emph{abstract neighborhood system} $\mathcal{V},$ i.e., a subset of positive elements that is directed downward,
  closed for addition and multiplication by  (strictly) positive scalars. The elements $v$ of $\mathcal{V}$ define \emph{upper} (\emph{lower}) \emph{neighborhoods} for the elements of
$\mathcal{P}$ by $ v(a)=\{b\in \mathcal{P}:b\leq a+v\}$
(respectively, $(a)v=\{b\in \mathcal{P}:a\leq b+v\}$),
 creating the \emph{upper}, respectively \emph{lower topologies} on $\mathcal{P}.$
  Their common refinement is called  the \emph{ symmetric topology.}
 We assume  all elements of $\mathcal{P}$ to be \textit{bounded
below}, i.e., for every $a\in \mathcal{P}$ and $v\in \mathcal{V}$ we
have $0\leq a+\rho v$ for some $\rho>0.$ Finally, a \textit{locally
convex cone} $(\mathcal{P},\mathcal{V})$ is a subcone of a full
locally convex cone, not necessarily containing the abstract
neighborhood system $\mathcal{V}.$ \par  For a locally convex cone
$(\mathcal{P},\mathcal{V)}$ the collection of all sets $
\widetilde{v}\subseteq\mathcal{P}^{2},$  where
$\widetilde{v}=\{(a,b):a\leq b+v\}$ for all  $v\in \mathcal{V},$
defines a \emph{convex quasi-uniform structure} on $\mathcal{P}.$ On
the other hand,  every convex quasi-uniform structure leads to a
full locally convex cone, including $\mathcal{P}$ as a subcone and
induces the same convex quasi-uniform structure. For details see
\cite [Ch I, 5.2]{KWR}.
\par  For cones $\mathcal{P}$ and $\mathcal{Q},$ a mapping  $t:\mathcal{P}\to\mathcal{Q}$ is called a \emph{linear operator},
if $t(a+b)=t(a)+t(b)$ and $t(\alpha a)=\alpha t(a)$ for all $a,b\in
\mathcal{P}$ and $\alpha\geq 0.$
 If $\mathcal{V}$ and
$\mathcal{W}$ are abstract neighborhood systems on $\mathcal{P}$
and $\mathcal{Q}$, a linear operator $t:\mathcal{P}\to\mathcal{Q}$
is called \emph{uniformly continuous}(\emph {u-continuous}), if
for every $w\in \mathcal{W}$ there is $v\in \mathcal{V}$ such that
$t(a)\leq t(b)+ w$ whenever $a\leq b+v$. Uniform continuity
implies continuity with respect to the upper, lower and symmetric
topologies on $\mathcal{P}$ and $\mathcal{Q}$. Endowed with the
neighborhood system
$\mathcal{V}=\{\varepsilon\in\mathbb{R}:\varepsilon>0\},$
$\overline{\mathbb{R}}$ is a full locally convex cone. The set of
all u-continuous linear functionals $\mu:\mathcal{P}\to
\overline{\mathbb{R}}$ is a cone called the \textit{dual cone} of
$\mathcal{P}$ and denoted by $\mathcal{P}^{\ast}.$ In a locally
convex cone $(\mathcal{P},\mathcal{V})$ the \textit{polar}
$v^{\circ}$ of $v\in\mathcal{V}$ is defined by
$v^{\circ}=\{\mu\in\mathcal{P}^{\ast}:a\leq b +v
\,\,\,\mbox{implies}\,\,\,\mu(a)\leq\mu(b)+1\}.$ Obviously we have
$\mathcal{P}^{\ast}=\cup_{v\in \mathcal{V}} v^{\circ}.$

\par A linear mapping $\Phi:\mathcal{P}\to\mathcal{Q}$ is
called an \emph{embedding }of $(\mathcal{P},\mathcal{V})$ into
$(\mathcal{Q},\mathcal{W})$ if it can be extended to a  mapping
$\Phi:\mathcal{P}\cup \mathcal{V}\to\mathcal{Q}\cup \mathcal{W}$
such that $\Phi(\mathcal{V})=\mathcal{W}$ and
$$
a\leq b+v\quad\mbox{holds if and only if
}\quad\Phi(a)\leq\Phi(b)+\Phi(v)
$$
for all $a,b\in \mathcal{P}$ and $v\in \mathcal{V}.$  This condition
implies that $\Phi$ is u-continuous, and in case that $\Phi$ is one
to one, the inverse operator $\Phi^{-1}:\Phi(\mathcal{P})\to
\mathcal{P}$ is also u-continuous. Embeddings are meant to preserve
not just the topological structure, but also the particular
neighborhood system of a locally convex cone \cite [Ch I, 2.2]{WRo}.
\section{Inductive limits}
\par Let
$\mathcal{P}_{\gamma},$ $\gamma\in\Gamma$ be  cones and put
$\mathcal{P}=\times_{\gamma\in\Gamma}\mathcal{P}_{\gamma}.$ For
elements $a,b\in\mathcal{P},$ $a=
\times_{\gamma\in\Gamma}a_{\gamma},b=\times_{\gamma\in\Gamma}b_{\gamma}$
and $\alpha\geq0$ we set
$a+b=\times_{\gamma\in\Gamma}\,(a_{\gamma}+b_{\gamma})$ and $\alpha
a=\times_{\gamma\in\Gamma}\,(\alpha a_{\gamma}).$ With these
operations
 $\mathcal{P}$ is a cone which is called the
 \emph{product cone} of $\mathcal{P}_{\gamma}.$ The subcone of the product cone
$\mathcal{P}$ spanned by $\cup \mathcal{P}_{\gamma}$ (more
precisely, by $\cup j_{\gamma}(\mathcal{P}_{\gamma})$, where
$j_{\gamma}:\mathcal{P}_{\gamma}\to\mathcal{P}$ is the injection
mapping) is said to be the \textit{direct sum cone} of
$\mathcal{P}_{\gamma}$ and denoted by
$\mathcal{Q}=\sum_{\gamma\in\Gamma}\mathcal{P}_{\gamma}.$ It is
worth remembering that here we only use positive scalars.
\par Suppose that
 $\texttt{Q}$  is a
 cone, for each $\gamma\in\Gamma,$
 $(\mathcal{P}_{\gamma},\leq_{\gamma})$ is an ordered cone and
$\mathrm{T}_{\gamma}:\mathcal{P}_{\gamma}\to\texttt{Q}$ is a linear
mapping such that
$\texttt{Q}=\sum_{\gamma\in\Gamma}\mathrm{T}_{\gamma}(\mathcal{P}_{\gamma}).$
For elements
$\texttt{a},\texttt{b}\in\texttt{Q},\texttt{a}=\sum_{\gamma\in\Delta}
\mathrm{T}_{\gamma}(a_{\gamma}),\texttt{b}=\sum_{\gamma\in\Theta}\mathrm{T}_{\gamma}(b_{\gamma}),$
$\Delta_{}=\{\gamma\in\Gamma:\mathrm{T}_{\gamma}(a_{\gamma})\neq
0\}$  and
$\Theta=\{\gamma\in\Gamma:\mathrm{T}_{\gamma}(b_{\gamma})\neq 0\},$
we set
 \begin{equation*}
 \texttt{a}\leq_{\overrightarrow{\Gamma}}
\texttt{b}
\end{equation*}\\
\emph{if for each $\gamma\in\Delta\cup\Theta$ there are
$c^{}_{\gamma_{1}},d_{\gamma_{1}},...,
c^{}_{\gamma_{n_{\gamma}}},d_{\gamma_{n_{\gamma}}}\in\mathcal{P}_{\gamma}$
 such that
$\mathrm{T}_{\gamma}(c_{\gamma_{1}})=\mathrm{T}_{\gamma}(a_{\gamma}),
\mathrm{T}_{\gamma}(d_{n_{\gamma}})=\mathrm{T}_{\gamma}(b),$  as
well as $c_{\gamma_{i}}\leq_{\gamma} d_{\gamma_{i}}$ for all
$i=1,2,...,n_{\gamma}$ and
$\mathrm{T}_{\gamma}(d_{\gamma_{i}})=\mathrm{T}_{\gamma}(c_{\gamma_{i+1}})$
for $ i=1,2,...,n_{\gamma}-1.$}

  \par As required, the order $\leq_{\overrightarrow{\Gamma}}$ is reflexive,
transitive and compatible with the algebraic operations of
$\texttt{Q}.$ This
 is the \emph{finest order} on
$\texttt{Q}$ such that the mappings
$a_{\gamma}\to\mathrm{T}_{\gamma}(a_{\gamma}):\mathcal{P}_{\gamma}\to\texttt{Q}
$  are monotone.
\par For each $\gamma\in\Gamma,$ let $(\mathcal{P}_{\gamma},\mathcal{V}_{\gamma})$ be a locally convex
cone and denote by
$\mathcal{W}=\times_{\gamma\in\Gamma}\mathcal{V}_{\gamma},$ the
product of the abstract neighborhoods $\mathcal{V}_{\gamma}.$ In
fact, $\mathcal{W}$ is a subcone (with out zero) of the product cone
of the all corresponding full cones containing
$\mathcal{V}_{\gamma}.$ The product
$\mathcal{W}=\times_{\gamma\in\Gamma}\mathcal{V}_{\gamma}$ leads us
to define the \emph{finest convex quasi-uniform structure}, as well
as, the \emph{finest locally convex cone topology} on $\texttt{Q}$
such that the  mappings $\mathrm{T}_{\gamma}$ are u-continuous:
\begin{definition}\label{def:2.1}
For  $\texttt{a},\texttt{b}\in\texttt{Q},$ $\texttt{a}=\sum_{\gamma\in\Delta}\mathrm{T}_{\gamma}(a_{\gamma}),$
$\texttt{b}=\sum_{\gamma\in\Theta}\mathrm{T}_{\gamma}(b_{\gamma})$
and $w\in \mathcal{W},w=\times_{\gamma\in\Gamma}v_{\gamma},$ we set
\begin{equation*}
\texttt{a}\leq_{\overrightarrow{\Gamma}} \texttt{b} +
\overrightarrow{w}
\end{equation*}
\emph{if for each $\gamma\in\Delta\cup\Theta$ there are
$c^{}_{\gamma_{1}},d_{\gamma_{1}},...,
c^{}_{\gamma_{n_{\gamma}}},d_{\gamma_{n_{\gamma}}}\in\mathcal{P}_{\gamma}$
 and $0\leq
\lambda_{\gamma_{1}},...,\lambda_{\gamma_{n_{\gamma}}} \in\mathbb{R}
$ such that $
\lambda_{\gamma_{1}}+...+\lambda_{\gamma_{n_{\gamma}}}\leq 1$  and
$\mathrm{T}_{\gamma}(c_{\gamma_{1}})=\mathrm{T}_{\gamma}(a_{\gamma}),
\mathrm{T}_{\gamma}(d_{\gamma_{n_{\gamma}}})=\mathrm{T}_{\gamma}(b),$
as well as $c_{\gamma_{i}}\leq_{\gamma}
d_{\gamma_{i}}+\lambda_{\gamma_{i}}v_{\gamma}$ for all
$i=1,2,...,n_{\gamma}$ and
$\mathrm{T}_{\gamma}(d_{\gamma_{i}})=\mathrm{T}_{\gamma}(c_{\gamma_{i+1}})$
for $ i=1,2,...,n_{\gamma}-1,$
$\sum_{\gamma\in\Delta\cup\Theta}\sum_{i=1}^{n_{\gamma}}\lambda^{}_{\gamma_{i}}\leq
1.$} \par  If, for each $w\in \mathcal{W},$ we put
$\widetilde{\overrightarrow{w}}=\{(\texttt{a},\texttt{b})\in
\texttt{Q}^{2}:\texttt{a}\leq_{\overrightarrow{\Gamma}}\texttt{b}+\overrightarrow{w}\},$
then
$\widetilde{\overrightarrow{\mathcal{W}}}=\{\widetilde{\overrightarrow{w}}:w\in\mathcal{W}\}$
forms the finest convex quasi-uniform structure on $\texttt{Q}$
which makes the mappings $\mathrm{T}_{\gamma},$ u-continuous. Then,
according to \cite [Ch I, 5.4]{KWR}, there exists a full cone
$\texttt{Q}\oplus\overrightarrow{\mathcal{W}_{0}},$ with abstract
neighborhood system
$\overrightarrow{\mathrm{W}}=\{0\}\oplus\overrightarrow{\mathcal{W}},$
whose neighborhoods yield the same convex quasi-uniform structure on
$\texttt{Q}.$ The elements $\overrightarrow{w}\in
\overrightarrow{\mathcal{W}},$
$w=\times_{\gamma\in\Gamma}v_{\gamma},$ form a basis for
$\overrightarrow{\mathrm{W}}$ in the following sense: For every
$\overrightarrow{\mathrm{w}}\in\overrightarrow{\mathrm{W}},$ there
is $w\in\mathcal{W}$ such that
$\texttt{a}\leq_{\overrightarrow{\Gamma}}
\texttt{b}+\overrightarrow{w}$ for $\texttt{a},\texttt{b}\in
\texttt{Q}$ implies that $\texttt{a}\leq_{\overrightarrow{\Gamma}}
\texttt{b}\oplus\overrightarrow{\mathrm{w}}.$ The locally convex
cone topology on $\texttt{Q}$ induced by
$\overrightarrow{\mathrm{W}}$ is called the \emph{locally convex
inductive limit cone of
$(\mathcal{P}_{\gamma},\mathcal{V}_{\gamma})$ by the linear mappings
$\mathrm{T}_{\gamma}$} and denoted by
$(\texttt{Q},\overrightarrow{\mathcal{W}})$ (cf. \cite [ 3.1]{WRl}).
\end{definition}
\begin{remark}
(i) If for each $\gamma\in\Gamma,$
$j_{\gamma}:\mathcal{P}_{\gamma}\to\mathcal{Q}$ is the injection
mapping, then
$\mathcal{Q}=\sum_{\gamma\in\Gamma}\mathcal{P}_{\gamma}$ will be the
locally convex inductive limit cone of
$(\mathcal{P}_{\gamma},\mathcal{V}_{\gamma})$ by  mappings
$j_{\gamma}$ which is called the\emph{ locally convex direct sum
cone of $(\mathcal{P}_{\gamma},\mathcal{V}_{\gamma})$} and denoted
by
$(\mathcal{Q},\mathcal{W})=(\sum_{\gamma\in\Gamma}\mathcal{P}_{\gamma},\times_{\gamma\in\Gamma}\mathcal{V}_{\gamma}).$
The direct sum neighborhood
$w\in\mathcal{W},w=\times_{\gamma\in\Gamma}v_{\gamma},$ for
$a,b\in\mathcal{Q},$
$a=\sum_{\gamma\in\Delta}a_{\gamma},b=\sum_{\gamma\in\Theta}b_{\gamma}$
is defined by
\begin{equation*}
a\leq_{\Gamma} b+w
\end{equation*}
\emph{if $a_{\gamma} \leq_{\gamma} b_{\gamma}+\alpha_{\gamma}
v_{\gamma}$ for all $\gamma\in\Delta\cup\Theta,$ where
$\sum_{\gamma\in\Delta\cup\Theta} \alpha_{\gamma}\leq 1$} (see \cite
[Definition 2.1]{MRlcpd}) .
 \par (ii) If for each $\gamma\in\Gamma,$ we
put
$\texttt{Q}_{\gamma}=\mathrm{T}_{\gamma}(\mathcal{P}_{\gamma}),$
then $\texttt{Q}$ will be the direct sum cone of
$\texttt{Q}_{\gamma},$ i.e.,
$\texttt{Q}=\sum_{\gamma\in\Gamma}\texttt{Q}_{\gamma}$ and
   $(\texttt{Q}_{\gamma},\overrightarrow{\mathcal{V}}_{\gamma})$ is the locally convex
inductive limit
  cone of $(\mathcal{P}_{\gamma},\mathcal{V}_{\gamma})$ by mapping
  $\mathrm{T}_{\gamma}.$ For elements $\texttt{a},\texttt{b}\in\texttt{Q},\texttt{a}=\sum_{\gamma\in\Delta}\texttt{a}_{\gamma},
  \texttt{b}=\sum_{\gamma\in\Theta}\texttt{b}_{\gamma}$ and $\overrightarrow{w}\in \overrightarrow{\mathcal{W}},$
   $w=\times_{\gamma\in\Gamma}v_{\gamma},$ we have
     $\texttt{a}\leq_{\overrightarrow{\Gamma}}\texttt{b}+\overrightarrow{w}$ if
   \begin{eqnarray*}
\texttt{a}_{\gamma}\leq_{\overrightarrow{\gamma}}
   \texttt{b}_{\gamma}+\sum_{i=1}^{n_{\gamma}}\lambda^{}_{\gamma_{i}}\overrightarrow{v}_{\gamma}\,\,\mbox{for
   all}\,\,\gamma\in\Delta\cup\Theta,\,\,\mbox{where}\,\,
\sum_{\gamma\in\Delta\cup\Theta}\sum_{i=1}^{n_{\gamma}}\lambda_{\gamma_{i}}\leq
1.
   \end{eqnarray*}
 That is,
$(\texttt{Q},\overrightarrow{\mathcal{W}})=(\sum_{\gamma\in\Gamma}
\texttt{Q}_{\gamma},\times_{\gamma\in\Gamma}\overrightarrow{\mathcal{V}}_{\gamma}),$
i.e., every inductive limit
$(\texttt{Q},\overrightarrow{\mathcal{W}})$ is the direct sum cone
of its components
$(\texttt{Q}_{\gamma},\overrightarrow{\mathcal{V}_{\gamma}}).$
\end{remark}
\par According to  \cite {WRl}, W. Roth  has
 introduced the notion of locally convex quotient cones as
 following: Suppose  $\sim$ is an equivalence relation on a locally convex
cone $(\mathcal{P},\mathcal{V})$ which is  \emph{compatible} with
the algebraic operations in $\mathcal{P},$ that is,  $a+c\sim b+c$
and $\alpha a\sim\alpha b$ whenever $a\sim b$ for $a,b,c\in
\mathcal{P}$ and $\alpha\geq0.$ The operations
$\widetilde{a}+\widetilde{b}=\widetilde{a+b}$ and $\alpha
\widetilde{a}=\widetilde{\alpha a}$ are well-defined for $a,b\in
\mathcal{P}$ and $\alpha\geq0,$ and $\widetilde{\mathcal{P}}
=\{\widetilde{a}:a\in\mathcal{P}\}$ becomes a cone with these
operations, where $\widetilde{a}$ is the \emph{equivalence class} of
$a.$ For $\widetilde{a},\widetilde{b}\in \widetilde{\mathcal{P}},$
and $v\in \mathcal{V}\cup\{0\}$ we set
\begin{eqnarray*}
\widetilde{a}\lesssim \widetilde{b} + \widetilde{v}
\end{eqnarray*}
\emph{if there are $c_{1},d_{1},..., c_{n},d_{n}\in\mathcal{P}$
 and $0\leq
\lambda_{1},...,\lambda_{n} \in\mathbb{R}$ such that $
\lambda_{1}+...+\lambda_{n}\leq 1$ and $c_{1}\sim a,$ $d_{n}\sim b,$
as well as $c_{i}\leq d_{i}+\lambda_{i}v$ for all $i=1,2,...,n$ and
$d_{i}\sim c_{i+1}$ for $ i=1,2,...,n-1,$ where
$\sum_{i=1}^{n}\lambda_{i}\leq 1.$}
 \par Then $\overset{\Gamma}{\sim}$
$\widetilde{\mathcal{V}}=\{\widetilde{v}:v\in\mathcal{V}\}$ defines
the finest locally convex cone topology on $\widetilde{\mathcal{P}}$
which makes the canonical projection
$\Pi(a)=\widetilde{a}:\mathcal{P}\to\widetilde{ \mathcal{P}},$
u-continuous. This is called the \emph{locally convex quotient cone
of $(\mathcal{P},\mathcal{V})$ over $\sim$} and denoted by
$(\widetilde{\mathcal{P}},\widetilde{\mathcal{V}}).$
\begin{remark} A quotient cone topology  $(\widetilde{\mathcal{P}},\widetilde{\mathcal{V}})$ is identical to   $(\Pi(\mathcal{P}),\overrightarrow{\mathcal{V}}),$
where $(\Pi(\mathcal{P}),\overrightarrow{\mathcal{V}})$ is the
locally convex inductive limit cone of $(\mathcal{P},\mathcal{V})$
by projection mapping $\Pi;$ for,  we have $\Pi(c_{1})=\Pi(a),$
$\Pi(d_{n})=\Pi(b)$ and $\Pi(d_{i})=\Pi(c_{i+1})$ for $i=1,2,...,n,$
hence $\widetilde{a}\lesssim \widetilde{b}+\widetilde{v}$ if and
only if $\Pi(a)\leq_{\overrightarrow{\Pi}}
\Pi(b)+\overrightarrow{v},$ where $\leq_{\overrightarrow{\Pi}}$ is
the inductive limit order on $\Pi(\mathcal{P}).$
\end{remark}
\begin{proposition} \label{prop:2.7} If  $(\texttt{Q},\overrightarrow{\mathcal{W}})$ is the locally convex
inductive limit cone of
$(\mathcal{P}_{\gamma},\mathcal{V}_{\gamma})$ by the linear mappings
$\mathrm{T}_{\gamma}$ and
$(\mathcal{Q},\mathcal{W})=(\sum_{\gamma\in\Gamma}\mathcal{P}_{\gamma},\times_{\gamma\in\Gamma}\mathcal{V}_{\gamma}),$
 then
\begin{enumerate}
\item [$(\mathrm{a})$] $(\texttt{Q},\overrightarrow{\mathcal{W}})$ is the locally convex inductive limit cone of
$(\mathcal{Q},\mathcal{W})$ by the linear mapping  $\mathrm{T}:
(\mathcal{Q},\mathcal{W})\to
(\texttt{Q},\overrightarrow{\mathcal{W}})$ such that
$\mathrm{T}(a)=\sum_{\gamma\in\Delta}\mathrm{T}_{\gamma}(a_{\gamma})$
for all $a\in\mathcal{Q},a=\sum_{\gamma\in\Delta}a_{\gamma};$ in
particular $\widetilde{\overrightarrow{\mathcal{W}}}$ is the finest
convex quasi-uniform structure on $\texttt{Q}$ such that the mapping
$\mathrm{T}$ is u-continuous,
\item [$(\mathrm{b})$] if $(\mathcal{R},\mathcal{Z})$ is a locally convex cone and $\Phi:(\texttt{Q},\overrightarrow{\mathcal{W}})\to
(\mathcal{R},\mathcal{Z})$ is a linear mapping, then $\Phi$ is
u-continuous if and only if
$\Phi\circ\mathrm{T}:(\mathcal{Q},\mathcal{W})\to
(\mathcal{R},\mathcal{Z})$ such that
$(\Phi\circ\mathrm{T})(a)=\sum_{\gamma\in\Delta}(\Phi\circ\mathrm{T}_{\gamma})(a_{\gamma})$
is u-continuous, that is,  for each $\gamma\in\Gamma,$
$\Phi\circ\mathrm{T}_{\gamma}
:(\mathcal{P}_{\gamma},\mathcal{V}_{\gamma})\to
(\mathcal{R},\mathcal{Z})$ is u-continuous.
 \end{enumerate}
\end{proposition}
\begin{proof} (a) According to  Definition \ref{def:2.1}, for
$\texttt{a},\texttt{b}\in\texttt{Q},$
$\texttt{a}=\sum_{\gamma\in\Delta}\mathrm{T}_{\gamma}(a_{\gamma}),$
$\texttt{b}=\sum_{\gamma\in\Theta}\mathrm{T}_{\gamma}(b_{\gamma})$
and $w\in \mathcal{W},w=\times_{\gamma\in\Gamma}v_{\gamma},$ we have
$\texttt{a}\leq_{\overrightarrow{\Gamma}}\texttt{b}+\overrightarrow{w}$
if and only if
 \begin{equation*}
 c_{i}\leq_{\Gamma}d_{i}+\lambda'_{i} w,\quad
 i=1,2,...,n_{(\Delta\cup\Theta)},\quad \lambda'_{i}=\sum_{\gamma\in\Delta\cup\Theta}\lambda_{\gamma_{i}}
\end{equation*}
\begin{equation*}
\mathrm{T}(c^{}_{1})=\mathrm{T}(a),\quad \mathrm{T}
(d_{n_{(\Delta\cup\Theta)}})=\mathrm{T}(b),\quad
\mathrm{T}(d^{}_{i})=\mathrm{T}(c^{}_{i+1}),\,
i=1,2,..,n_{\Delta\cup\Theta}-1,
 \end{equation*}
 where $\sum_{i=1}^{n(\Delta\cup\Theta)}\lambda '_{i}\leq 1$ and  $c_{1},d_{1},...,$ $c_{n_{(\Delta\cup\Theta)}},\,d_{n_{(\Delta\cup\Theta)}}\in
 \mathcal{Q},$
\begin{equation}\label{equ:(2.1)}
c_{1}=\sum_{\gamma\in\Delta\cup\Theta}c_{\gamma_{1}},\,\,
d_{i}=\sum_{\gamma\in\Delta\cup\Theta} d_{\gamma_{i}},\,\,
c_{i+1}=\sum_{\gamma\in\Delta\cup\Theta}c_{\gamma_{i+1}},\,\,
d_{n_{(\Delta\cup\Theta)}}=\sum_{\gamma\in\Delta\cup\Theta}
d_{n_{\gamma}}.
\end{equation}
\par  For (b), if $\Phi$ is u-continuous, then $\Phi\circ\mathrm{T}$
is u-continuous by (a). Conversely, let for each $z\in\mathcal{Z},$
there is  $w\in \mathcal{W},w=\times_{\gamma\in\Gamma}v_{\gamma}$
such that $a\leq_{\Gamma} b+ w$ implies
$(\Phi\circ\mathrm{T})(a)\leq (\Phi\circ\mathrm{T})(b)+z.$
 For $\texttt{a},\texttt{b}\in \texttt{Q},$ $\texttt{a}=\sum_{\gamma\in\Delta}\mathrm{T}_{\gamma}(a_{\gamma_{}}),$
$\texttt{b}=\sum_{\gamma\in\Theta}\mathrm{T}_{\gamma}(b_{\gamma_{}}),$
let $\texttt{a}\leq_{\overrightarrow{\Gamma}}\texttt{b}+
\overrightarrow{w}.$ If we choose the elements
$c_{1},d_{1},...,c_{n_{(\Delta\cup\Theta)}},d_{n_{(\Delta\cup\Theta)}}\in\mathcal{Q}$
as  (\ref{equ:(2.1)}), then
$$
(\Phi\circ\mathrm{T})(c_{i})\leq_{}(\Phi\circ\mathrm{T})(d_{i})+\lambda'_{i}
\,z\quad(i=1,2,...,n_{(\Delta\cup\Theta)})
$$
which yields
$$\Phi(\texttt{a})=(\Phi\circ\mathrm{T})(c_{1})\leq
(\Phi\circ\mathrm{T})(d_{n_{(\Delta\cup\Theta)}})
+\sum_{\gamma\in\Delta\cup\Theta}
\lambda'_{i}z\leq\Phi(\texttt{b})+z,$$ i.e., $\Phi$ is u-continuous.
\end{proof}
\begin{corollary}\label{Cor:2.10}
If
$\overrightarrow{w}\in\overrightarrow{\mathcal{W}},w=\times_{\gamma\in\Gamma}v_{\gamma},$
then \begin{enumerate}
\item [$(\mathrm{a})$]
$\mu\in{\overrightarrow{w}}^{\circ}$ if and only if $\mu\circ
\mathrm{T}\in w^{\circ},$ i.e.,  $\mu\circ \mathrm{T}_{\gamma}\in
v_{\gamma}^{\circ}\,\,\,\,\mbox{for all}\,\,\gamma\in\Gamma,$
\item [$(\mathrm{b})$]
${\overrightarrow{w}}^{\circ}=\times_{\gamma\in\Gamma}{\overrightarrow{v}}^{\circ}_{\gamma};$
in particular,
$\texttt{Q}^{\ast}=\times_{\gamma\in\Gamma}\texttt{Q}^{\ast}_{\gamma}.$
\end{enumerate}
\end{corollary}
\begin{proof}
Part (a) holds by Proposition \ref{prop:2.7} (b). For (b), let
$\mu\in {\overrightarrow{w}}^{\circ},$
$\mu=\times_{\gamma\in\Gamma}\mu_{\gamma}.$ By  (a), for each
$\overrightarrow{v}_{\gamma}\in\overrightarrow{\mathcal{V}}_{\gamma},$
we have $\mu_{\gamma}\in {\overrightarrow{v}}^{\circ}_{\gamma}$ if
and only if
$\mu\circ\mathrm{T}_{\gamma}=\mu_{\gamma}\circ\mathrm{T}_{\gamma}\in
v^{\circ}_{\gamma}.$ That is,  $ \mu\in
{\overrightarrow{w}}^{\circ}$ if and only if $\mu_{\gamma}\in
{\overrightarrow{v}}^{\circ}_{\gamma}$ for all $\gamma\in\Gamma.$
\end{proof}
\par  A locally convex cone $(\mathcal{P},\mathcal{V})$ is said to
have the  \textit{strict separation property}, in short (SP), if for
all $a,b\in \mathcal{P}$ and $v\in \mathcal{V}$ with $a\not\leq b +
\rho v$ for some $\rho> 1,$ there is a $\mu\in v^{\circ}$ such that
$\mu(a)>\mu(b)+1.$
\begin{corollary} \label{Cor:2.11}
\par    If  the all
$(\mathcal{P}_{\gamma},\mathcal{V}_{\gamma})$ have  the strict
separation  property, then
$(\texttt{Q},\overrightarrow{\mathcal{W}})$ will also carry the same
property.
\end{corollary}
\begin{proof} Let $\texttt{a},\texttt{b}\in\texttt{Q}$ and
$\overrightarrow{w}\in\overrightarrow{\mathcal{W}},w=\times_{\gamma\in\Gamma}v_{\gamma}$
such that $\texttt{a}\not\leq_{\overrightarrow{\Gamma}}
\texttt{b}+ \rho\, \overrightarrow{w}$ for some $\rho>1.$ Then
$a_{\gamma_{}}\not\leq_{\gamma} b_{\gamma_{}}+\rho\,
v_{\gamma_{}}$ for some $\gamma\in\Delta\cup \Theta,$ hence by the
strict separation  property of
$(\mathcal{P}_{\gamma},\mathcal{V}_{\gamma}),$ there is
$\mu_{\gamma_{}}\in v^{\circ}_{\gamma_{}}$ such that
$\mu_{\gamma_{}}(a_{\gamma_{}})>\mu_{\gamma_{}}
(b_{\gamma_{}})+1.$ Define the mapping
$\nu_{\gamma}:\texttt{Q}_{\gamma}\to\overline{\mathbb{R}}$  by
$\nu_{\gamma}(\texttt{a}_{\gamma})=\mu_{\gamma}(a_{\gamma})$ for
all $\texttt{a}_{\gamma}\in\texttt{Q}_{\gamma}.$ Obviously, we
have $\nu_{\gamma}\in {\overrightarrow{v}}^{\circ}_{\gamma}.$ If
we put $\nu=\times_{\lambda\in\Gamma}\nu_{\lambda},$ where
$\nu_{\lambda}=\nu_{\gamma_{}}$ for $\lambda=\gamma$ and
$\nu_{\lambda}=0$ otherwise, then
$\nu(\texttt{a})>\nu(\texttt{b})+1$  and $\nu\in
{\overrightarrow{w}}^{\circ}$ by Corollary \ref{Cor:2.10} (b).
\end{proof}
\begin{remark} \label{Rem:2.12} If $(\mathcal{Q},\mathcal{W})=(\sum_{\gamma\in\Gamma}
\mathcal{P}_{\gamma},\times_{\gamma\in\Gamma}\mathcal{V}_{\gamma})$
then, by Corollary \ref{Cor:2.10} (b), for each
$w\in\mathcal{W},w=\times_{\gamma\in\Gamma}v_{\gamma},$ we have
$w^{\circ}=\times_{\gamma\in\Gamma}v^{\circ}_{\gamma};$ in
particular,
$\mathcal{Q}^{\ast}=\times_{\gamma\in\Gamma}\mathcal{P}^{\ast}_{\gamma}.$
Also, according to  Corollary \ref{Cor:2.11},
 $(\mathcal{Q},\mathcal{W})$ will carry the strict separation property, whenever its components
$(\mathcal{P}_{\gamma},\mathcal{V}_{\gamma})$ have  the same
property.
\end{remark}
\par  Suppose  $(\mathcal{P}, \mathcal{V})$ is  a full locally convex
cone and let $Conv(\mathcal{P})$ be the cone of all non-empty convex
subsets of $\mathcal{P}$  with  usual addition and scalar
multiplication of sets by non-negative scalars. If we identify the
elements of $\mathcal{V}$ with singleton sets $\bar{v}=\{v\},$ then
$\overline{\mathcal{V}}=\{\bar{v}:v\in\mathcal{V}\}$ is a subset of
$Conv(\mathcal{P}),$ which can be preordered  using the preorder of
$\mathcal{P}.$ For $A,B\in Conv(\mathcal{P}),$ we define $A\leq B$
if for every $a\in A$ there is $b\in B$ with  $a\leq b.$ Since its
elements are bounded below  as are the elements of $\mathcal{P},$
$(Conv(\mathcal{P}),\overline{\mathcal{V}})$ becomes a full locally
convex cone.  If  $\overline{Conv}(\mathcal{P})$ is the set of all
non-empty lower closed convex subsets of $\mathcal{P},$ then with
the following modified addition it will become a cone as well:
\begin{equation*}
A\overline{\oplus}B=\overline{A+B}\quad \mbox{for}\,\,\,A,B\in
\overline{Conv}(\mathcal{P}),
\end{equation*}
and $(\overline{Conv}(\mathcal{P}),$ $\overline{\mathcal{V}})$ forms
a locally convex cone. In particular,
$(\overline{\mathcal{P}},\overline{\mathcal{V}})$ is a locally
convex cone, where
$\overline{\mathcal{P}}=\{\overline{a}:a\in\mathcal{P}\}.$ For
details see \cite {KWR}.
\begin{remark} \label{Rem:2.4} (i) Let $(E,\mathcal{V},\leq)$ be a
locally convex ordered topological vector space, where $\mathcal{V}$
is a basis of closed, convex, balanced and order convex
neighborhoods of the origin in $E.$ The order on $Conv(E)$ is given
by $A\leq B$ if $A_{}\subset B_{}+E_{-},$ where $E_{-}=\{x\in
E:x\leq 0\}$ is the negative cone in $E.$ For every $A\in Conv(E)$
and $V\in\mathcal{V}$  there is a $\rho>0$ such that $\rho V\cap
A\neq\emptyset,$ that is, $0\in A+\rho V$ or $0\leq A+\rho V.$ Thus
every element of $Conv(E)$ is bonded below and
$(Conv(E),\mathcal{V})$ is a full locally convex cone. The embedding
$x\to \{x\}: E\to Conv(E)$ preserves the order structure of $E$ and
on it's image in $Conv(E)$ the three (upper, lower and symmetric)
topologies coincide with the given topology on $E.$ Thus
$(E,\mathcal{V})$ is a locally convex cone, but not a full cone (see
\cite [ 2.1 (c)] {WRl}. Now, for each $\gamma\in\Gamma,$ let
$(E_{\gamma},\mathcal{V}_{\gamma})$ be  a locally convex topological
vector space, $\texttt{E}$   a vector space and let
$\mathrm{T}_{\gamma}:(E_{\gamma},\mathcal{V}_{\gamma})\to
\texttt{E}$ be a linear mapping such that
$\texttt{E}=\sum_{\gamma\in\Gamma}\mathrm{T}_{\gamma}(E_{\gamma}).$
Suppose $(\texttt{E},\overrightarrow{\mathcal{W}})$ is the locally
convex inductive limit cone of $(E_{\gamma},\mathcal{V}_{\gamma})$
by $\mathrm{T}_{\gamma}$ and let $\mathcal{U}$ be a basis of closed
convex neighborhoods for inductive limit topology on $\texttt{E}.$
 For $\texttt{a},\texttt{b}\in\texttt{E},$
$\texttt{a}=\sum_{\gamma\in\Delta}\mathrm{T}_{\gamma}(a_{\gamma}),\texttt{b}=\sum_{\gamma\in\Theta}\mathrm{T}_{\gamma}(b_{\gamma})$
and neighborhood $\mathrm{U}\in\mathcal{U},$ we have
$\texttt{a}\in\texttt{b}+\mathrm{U}$ if and only if
\begin{eqnarray*}
a_{\gamma}\leq_{\gamma}
b_{\gamma}+\alpha_{\gamma}\mathrm{T}^{-1}_{\gamma}(\mathrm{U})\quad
\mbox{and}\quad b_{\gamma}\leq_{\gamma}
a_{\gamma}+\alpha_{\gamma}\mathrm{T}^{-1}_{\gamma}(\mathrm{U})
\end{eqnarray*}
 for all $\gamma\in\Delta\cup\Theta,$ where $\sum_{\gamma\in\Delta\cup\Theta}\alpha_{\gamma}\leq 1,$
 that is,
 $\texttt{a}\leq_{\overrightarrow{\Gamma}} \texttt{b}+\overrightarrow{w}$ and  $\texttt{b}\leq_{\overrightarrow{\Gamma}}
 \texttt{a}+\overrightarrow{w},$  where $w=\times_{\gamma\in\Gamma}\mathrm{T}^{-1}_{\gamma}(\mathrm{U})\in\mathcal{W}$
 or
$\mathrm{U}\subseteq
\overrightarrow{w}(\texttt{b})\overrightarrow{w}.$ Thus,
$\mathcal{U}$ is finer than the symmetric topology of
$(\texttt{E},\overrightarrow{\mathcal{W}}).$ Conversely, let
$w\in\mathcal{W},w=\times_{\gamma\in\Gamma}V_{\gamma}$ and let
$\mathrm{U}$ be the closed convex hull of
$\times_{\gamma\in\Gamma}\mathrm{T}_{\gamma}(V_{\gamma}).$ If
$\texttt{a}\leq_{\overrightarrow{\Gamma}}\texttt{b}+\overrightarrow{w},$
then $\texttt{a}\in\texttt{b}+\mathrm{U},$ so
$\overrightarrow{w}(\texttt{b})\overrightarrow{w}\subseteq
\mathrm{U}.$ Hence,  the symmetric topology of
$(\texttt{E},\overrightarrow{\mathcal{W}})$ is identical to the
inductive limit topology of
 $(E_{\gamma},\mathcal{V}_{\gamma})$ in the sense of locally convex topological
 vector spaces.
\par (ii) Let $(\texttt{Q},\overrightarrow{\mathcal{W}})$ be
the locally convex inductive limit cone of
$(\mathcal{P}_{\gamma},\mathcal{V}_{\gamma})$ by the linear mappings
$\mathrm{T}_{\gamma}$ and
$(\mathcal{Q},\mathcal{W})=(\sum_{\gamma\in\Gamma}\mathcal{P}_{\gamma},\times_{\gamma\in\Gamma}\mathcal{V}_{\gamma}).$
 We extend each $\mathrm{T}_{\gamma}$ to
$Conv(\mathcal{P}_{\gamma})$ by defining
$\overline{\mathrm{T}}_{\gamma}(A_{\gamma})=\overline{\mathrm{T}_{\gamma}(A_{\gamma})}$
 and obtain the linear mappings
$\overline{\mathrm{T}}_{\gamma}:Conv(\mathcal{P}_{\gamma})\to
\overline{Conv}(\texttt{Q}_{\gamma}).$
 If $\texttt{A}\in
\sum_{\gamma\in\Gamma}Conv(\texttt{Q}_{\gamma}),$ then
$\texttt{A}=\mathrm{T}(A)=\sum_{\gamma\in\Delta}\mathrm{T}_{\gamma}(A_{\gamma})$
for some $A\subseteq\mathcal{Q}, A=\sum_{\gamma\in\Delta}A_{\gamma}$
so by \cite [Proposition 3.2] {MRHpd}, we have
$\overline{\texttt{A}}=\sum_{\gamma\in\Delta}\overline{\mathrm{T}_{\gamma}(A_{\gamma})},$
hence $\sum_{\gamma\in\Gamma}\overline{Conv}(\texttt{Q}_{\gamma})
=\sum_{\gamma\in\Gamma}\overline{\mathrm{T}}_{\gamma}(Conv(\mathcal{P}_{\gamma}).$
For elements $\texttt{A},\texttt{B}\in
\sum_{\gamma\in\Gamma}Conv(\texttt{Q}_{\gamma}),$
$\texttt{A}=\sum_{\gamma\in\Delta}\mathrm{T}_{\gamma}(A_{\gamma}),$
$\,\texttt{B}=\sum_{\gamma\in\Theta}\mathrm{T}_{\gamma}(B_{\gamma})$
and $w\in\mathcal{W},w=\times_{\gamma\in\Gamma}v_{\gamma},$ we have
$\overline{\texttt{A}}\leq
\overline{\texttt{B}}+\overline{\overrightarrow{w}}$ if and only if
\begin{equation*}
\sum_{\gamma\in\Delta}\overline{\mathrm{T}}_{\gamma}(A_{\gamma})\leq_{\overrightarrow{\Gamma}}
\sum_{\gamma\in\Theta}\overline{\mathrm{T}}_{\gamma}(B_{\gamma})+\overrightarrow{\times}_
{\gamma\in\Gamma}\overline{v}_{\gamma},
\end{equation*}
 that is, $\overline{\overrightarrow{\mathcal{W}}}$
is identical to  neighborhood system
$\overrightarrow{\times}_{\gamma\in\Gamma}\overline{\mathcal{V}}_{\gamma},$
i.e.,
$(\sum_{\gamma\in\Gamma}\overline{Conv}(\texttt{Q}_{\gamma}),\overline{\overrightarrow{\mathcal{W}}})$
carries the inductive limit cone of
$(Conv(\mathcal{P}_{\gamma}),\overline{\mathcal{V}}_{\gamma})$ by
mappings $\overline{\mathrm{T}}_{\gamma}.$ In particular, for each
$a\in\mathcal{Q},a=\sum_{\gamma\in\Delta}a_{\gamma}$  we have
$\overline{\texttt{a}}=\sum_{\gamma\in\Delta}\overline{\texttt{a}}_{\gamma}$
 by \cite [Corollary
 3.3] {MRHpd} which
yields
$\overline{\texttt{Q}}=\sum_{\gamma\in\Gamma}\overline{\texttt{Q}}_{\gamma},$
so
$(\overline{\texttt{Q}},\overline{\overrightarrow{\mathcal{W}}})$
is the locally convex inductive limit cone of
$(\mathcal{P}_{\gamma},\mathcal{V}_{\gamma})$ by linear mappings
$\overline{\mathrm{T}}_{\gamma}:\mathcal{P}_{\gamma}\to
\overline{\texttt{Q}}$ such that
$\overline{\mathrm{T}}_{\gamma}(a_{\gamma})=\overline{\mathrm{T}_{\gamma}(a_{\gamma})}$
for all $a_{\gamma}\in\mathcal{P}_{\gamma}.$
\par (iii)
The locally convex cone topology
$(\overline{Conv}(\texttt{Q}),\overline{\overrightarrow{\mathcal{W}}})$
carries the  inductive limit of
$(Conv(\mathcal{Q}),\overline{\mathcal{W}})$ by linear mapping
$A\to\overline{\mathrm{T}}(A):Conv(\mathcal{Q})\to\overline{Conv}(\texttt{Q}),$
where $\overline{\mathrm{T}}(A)=\overline{\mathrm{T}(A)}.$ For, by
Proposition \ref{prop:2.7} (a) and \cite[Ch II, 1.6 (e)]{KWR}, the
mapping $\overline{\mathrm{T}}$ is u-continuous, so
$\overrightarrow{\overline{\mathcal{W}}}$ is finer that
$\overline{\overrightarrow{\mathcal{W}}}.$ On the other hand, for
elements $A,B\in Conv(\mathcal{Q})$ and $w\in\mathcal{W},$ if
$\overline{\mathrm{T}}(A)\leq\overline{\mathrm{T}}
(B)+\overline{\overrightarrow{w}},$ then for each $a\in A,$ there
is $b\in B$ such that
$\mathrm{T}(a)\leq_{\overrightarrow{\Gamma}}\mathrm{T}(b)+\overrightarrow{w},$
i.e.,
$\overline{\mathrm{T}}(A)\leq_{\overrightarrow{\overline{\mathrm{T}}}}\overline{\mathrm{T}}(B)+\overrightarrow{w},$
where $\leq_{\overrightarrow{\overline{\mathrm{T}}}}$ is the
inductive limit order on $\overline{Conv}(\texttt{Q}),$ so
$\overline{\overrightarrow{\mathcal{W}}}$ is identical to
$\overrightarrow{\overline{\mathcal{W}}}.$
 In particular,
 $(\overline{\texttt{Q}},\overline{\overrightarrow{\mathcal{W}}})$
is  the locally convex inductive limit cone of
 $(\mathcal{Q},\mathcal{W})$
 by linear mapping $\overline{\mathrm{T}}:\mathcal{Q}\to \overline{\texttt{Q}};$
 where $\overline{\mathrm{T}}(a)
 =\overline{\mathrm{T}(a)}$  for all $a\in\mathcal{Q}.$
\end{remark}
 \begin{example} \label{Exam:2.5} \par (i) If we consider the neighborhood system
 $\{\varepsilon>0:\varepsilon\in\mathbb{R}\}$ on
 $\overline{\mathbb{R}},$ then
 for each $n\in\mathbb{N},$ the direct sum neighborhood system on
$\overline{\mathbb{R}}^{n}$
 is given by  $\mathcal{W}_{n}=\{\varepsilon\,
1_{\mathbb{R}^{n}}:\varepsilon>0\}$ for all $n\in\mathbb{N},$
 where
$1_{\mathbb{R}^{n}}=(1,1,...,1).$ For elements
$a_{n},b_{n}\in\overline{\mathbb{R}}^{n},$
$a_{n}=(a^{1}_{n},...,a^{n}_{n}),\,b_{n}=(b^{1}_{n},...,b^{n}_{n}),$
 the neighborhood $w_{n}=\varepsilon \,1_{\mathbb{R}^{n}}$ for $\varepsilon>0,$
is defined by $a_{n}\leq_{n} b_{n}+\varepsilon \,1_{\mathbb{R}^{n}}$
if $a^{i}_{n}\leq b^{i}_{n}+\alpha _{i}\,\varepsilon,$
$i=1,2,...,n;$ where $\sum_{i=1}^{n}\alpha_{i}\leq 1.$ Let
$\mathcal{W}=\times_{n\in\mathbb{N}}\mathcal{W}_{n}$ be the
neighborhood system on  direct sum cone
$\sum_{n\in\mathbb{N}}\overline{\mathbb{R}}^{n}.$
 If for each
$n\in\mathbb{N},$ we define the linear mapping $\mathrm{T}_{n}$ on
$\overline{\mathbb{R}}^{n}$ for all
$a_{n}\in\overline{\mathbb{R}}^{n}$ by
$\mathrm{T}_{n}(a_{n})=\sum_{i=1}^{n}a^{i}_{n}\in\overline{\mathbb{R}},$
then
$\sum_{n\in\mathbb{N}}\mathrm{T}_{n}(\overline{\mathbb{R}}^{n})=\sum_{n\in\mathbb{N}}\overline{\mathbb{R}}:=\texttt{Q}.$
Thus, for elements $\texttt{a},\texttt{b}\in\texttt{Q},$
$\texttt{a}=\sum_{n\in\Delta}\mathrm{T}_{n_{}}(a_{n_{}})=\sum_{n\in\Delta}\sum_{i=1}^{n_{}}a^{i}_{n_{}},$
$\texttt{b}=\sum_{n\in\Theta}\mathrm{T}_{n_{}}(b_{n_{}})
=\sum_{n\in\Theta}\sum_{i=1}^{n_{}}b^{i}_{n_{}},$ and
$w\in\mathcal{W},$ $w=\varepsilon
\times_{n\in\mathbb{N}}1_{\mathbb{R}^{n}},$  $\varepsilon>0,$ the
inductive limit neighborhood on $\texttt{Q}$ is defined by
$\texttt{a}\leq_{\overrightarrow{\mathbb{N}}}
\texttt{b}+\varepsilon\,
\overrightarrow{\times}_{n\in\mathbb{N}}1_{\mathbb{R}^{n}}$ if
\begin{equation*}
 \sum_{i=1}^{n_{}}a^{i}_{n_{}}\leq
\sum_{i=1}^{n_{}}b^{i}_{n_{}}+\varepsilon/\sum_{n^{}\in\Delta\cup\Theta}\,n\quad\mbox{for
all}\,\,\,\,n\in\Delta\cup\Theta,
\end{equation*}
where $a^{i}_{n_{}}=b^{i}_{n_{}}=0$ for all
$n\in(\Delta\setminus\Theta)\cup(\Theta\setminus \Delta)$ and
$i=1,2,...,n.$ That is,
$(\sum_{n\in\mathbb{N}}\overline{\mathbb{R}},\overrightarrow{\mathcal{W}})$
is the locally convex inductive limit cone of
$(\overline{\mathbb{R}}^{n},\mathcal{W}_{n})$ by mappings
$\mathrm{T}_{n},$ where
$\overrightarrow{\mathcal{W}}=\{\varepsilon\,\overrightarrow{\times}_{n\in\mathbb{N}}1_{\mathbb{R}^{n}}:\varepsilon>0\}.$
\par (ii)
 According to \cite[Corollary 3.3]{MRHpd}, we have
 $\overline{\sum}_{n\in\mathbb{N}}\overline{\mathbb{R}}=\sum_{n\in\mathbb{N}}\overline{\overline{\mathbb{R}}},$
 where $\overline{\overline{\mathbb{R}}}=\{\overline{a}:a\in\overline{\mathbb{R}}\}$ and $\overline{a}=(-\infty,a]$ for all $a\in\overline{\mathbb{R}}.$
 We note that $\overline{Conv}(\overline{\mathbb{R}})$ is identical to
 $\overline{\overline{\mathbb{R}}},$ so
for all  $A_{n}\in Conv(\overline{\mathbb{R}}^{n}),$ we have
$\overline{\mathrm{T}}_{n}(A_{n})=\overline{\mathrm{T}_{n}(A_{n})}\in\overline{\overline{\mathbb{R}}},$
 whence
$\sum_{n\in\mathbb{N}}\overline{\mathrm{T}}_{n}(Conv(\overline{\mathbb{R}}^{n}))=\overline{\sum}_{n\in\mathbb{N}}\overline{\mathbb{R}}:=\texttt{Q}.$
Suppose
$\mathcal{W}=\times_{n\in\mathbb{N}}\overline{\mathcal{W}}_{n}$ is
the direct sum neighborhood system on
$\sum_{n\in\mathbb{N}}Conv(\overline{\mathbb{R}}^{n}),$ where
$\overline{\mathcal{W}}_{n}=\{\varepsilon\,
\overline{1}_{\mathbb{R}_{n}}:\varepsilon>0\}$ is the corresponding
neighborhood system on $Conv(\overline{\mathbb{R}}^{n})$ for all
$n\in\mathbb{N}.$  For elements
$\texttt{a},\texttt{b}\in\texttt{Q},$
$\texttt{a}=\sum_{n\in\Delta}\overline{\mathrm{T}}_{n_{}}(A_{n_{}}),$
$\texttt{b}=\sum_{n\in\Theta}\overline{\mathrm{T}}_{n_{}}(B_{n_{}}),
$ there exist $a_{n},$ $b_{n}\in\overline{\mathbb{R}}^{n}$ such
$\mathrm{T}_{n}(A_{n})=\mathrm{T}_{n}(a_{n})$ for all $n\in\Delta$
and $\mathrm{T}_{n}(B_{n})=\mathrm{T}_{n}(b_{n})$ for all
$n\in\Theta,$ so the inductive limit neighborhood on $\texttt{Q}$
for $\varepsilon>0$ is given  by
$\texttt{a}\leq_{\overrightarrow{\mathbb{N}}}
\texttt{b}+\varepsilon\,\overrightarrow{\times}_{n\in\mathbb{N}}\overline{1}_{\mathbb{R}_{n}}$
if
\begin{equation*}\overline{\sum}_{i=1}^{n}
a^{i}_{n}\subseteq \overline{\sum}_{i=1}^{n}
b^{i}_{n}+\mathbb{R}_{-}+\varepsilon/\sum_{n\in\Delta\cup\Theta}\,n
\end{equation*}
for all $n\in\Delta\cup\Theta.$ 
Thus
$(\overline{\sum}_{n\in\mathbb{N}}\overline{\mathbb{R}},\overrightarrow{\mathcal{W}})$
is the locally convex inductive limit cone of
$(Conv(\overline{\mathbb{R}}^{n}),\overline{\mathcal{W}}_{n})$ by
mappings $\overline{\mathrm{T}}_{n},$ where
$\overrightarrow{\mathcal{W}}=\{\varepsilon\,\overrightarrow{\times}_{n\in\mathbb{N}}\overline{1}_{\mathbb{R}^{n}}:\varepsilon>0\}.$
\par (iii) By  (i) and Remark 2.8 (iii), we infer that
$(\overline{Conv}(\sum_{n\in\mathbb{N}}\overline{\mathbb{R}}),\overline{\overrightarrow{\mathcal{W}}})$
carries the  inductive limit  of $(Conv(
\overline{\mathbb{R}}^{n}),\overline{\mathcal{W}})$ by mapping
$\overline{\mathrm{T}},$  where
$\overline{\mathrm{T}}(A)=\overline{\mathrm{T}(A)
}\in\overline{Conv}(\sum_{n\in\mathbb{N}}\overline{\mathbb{R} })$
 for all $A\in Conv(
\overline{\mathbb{R}}^{n}).$ For $w\in\mathcal{W},$
$w=\varepsilon\times_{n\in\mathbb{N}}1_{\mathbb{R}^{n}},$ we have
$\overline{\mathrm{T}}(A)\leq_{\overrightarrow{\overline{\mathrm{T}}}}\overline{\mathrm{T}}(B)+\overrightarrow{w}$
if $\texttt{A}\subseteq\texttt{
B}+\mathbb{R}^{\mathbb{N}}_{-}+\varepsilon
\times_{n\in\mathbb{N}}\mathbb{B}_{n},$ where $\mathbb{B}_{n}$ is
the closed unit ball of $\mathbb{R}^{n}$ for all $n\in\mathbb{N}.$
In particular,
$(\overline{\sum}_{n\in\mathbb{N}}\overline{\mathbb{R}},\overline{\overrightarrow{\mathcal{W}}})$
is the inductive limit of $(\sum_{n\in\mathbb{N}}
\overline{\mathbb{R}}^{n},\overline{\mathcal{W}})$ by mapping
$\overline{\mathrm{T}},$ where
$\overline{\mathrm{T}}(a)=\sum_{n\in\Delta}\overline{\sum}_{i=1}^{n_{}}a^{i}_{n_{}}$
for all $a\in\overline{\mathbb{R}}^{n},$
$a=\sum_{n\in\Delta}a_{n}.$
\end{example}
 \par  Let  $X$ be a
  set and let
$\mathcal{F}(X,\overline{\mathbb{R}})$ be the cone of all
$\overline{\mathbb{R}}$-valued functions on $X$ endowed with the
pointwise operations and order. A neighborhood system
$\varepsilon_{\mathcal{Y}}$ for
$\mathcal{F}(X,\overline{\mathbb{R}})$ may be created using a
suitable family $\mathcal{Y}$ of subsets in $X,$ directed downward
with respect to set inclusion, and  the neighborhoods
$\varepsilon_{Y}$ for $\varepsilon>0$ and $Y\in\mathcal{Y},$ defined
for functions $f,g\in\mathcal{F}(X,\overline{\mathbb{R}})$ as
 $f\leq g+\varepsilon_{Y}$
 if $f(x)\leq g(x)+\varepsilon$ for all $x\in Y.$ In this case we
 consider the subcone $\mathcal{F}_{b_{\mathcal{Y}}}(X,\overline{\mathbb{R}})$
 of all functions  $f$ in $\mathcal{F}(X,\overline{\mathbb{R}})$ that
  are
uniformly bounded below on the sets in $\mathcal{Y},$ i.e., for
every $\varepsilon>0$ and $Y\in\mathcal{Y},$ there is a $\rho>0$
such that $0\leq f+\rho \varepsilon_{Y}.$ Together with the
neighborhood system $\varepsilon_{\mathcal{Y}},$ it forms a locally
convex cone. For details see \cite [3.7 (f)]{WRl}.
 \begin{example}  For each $\gamma\in\Gamma,$ let  $X_{\gamma}$ be
 a set and  $\mathcal{Y}_{\gamma}$ a family of subsets of
 $X_{\gamma},$ directed downward with respect to set inclusion and
 let
$\mathcal{P}_{\gamma}=\mathcal{F}_{b_{\mathcal{Y}_{\gamma}}}(X_{\gamma},\overline{\mathbb{R}})$
 endowed with the neighborhood system $\varepsilon_{\mathcal{Y}_{\gamma}}.$
If for each $f_{\gamma}\in\mathcal{P}_{\gamma},$ we set
 $\mathcal{I}_{f_{\gamma}}=\{x_{\gamma}\in X_{\gamma}:f_{\gamma}(x_{\gamma})=+\infty\}$  and
 put $\texttt{Q}_{\gamma}=\{\mathcal{I}_{f_{\gamma}}:f_{\gamma}\in\mathcal{P}_{\gamma}\}$ for all $\gamma\in\Gamma,$ then
 with the
following operations $\texttt{Q}_{\gamma}$ becomes a cone:
\begin{equation*}
\mathcal{I}_{f_{\gamma}}\oplus
\mathcal{I}_{g_{\gamma}}=\mathcal{I}_{f_{\gamma}}\cup
\mathcal{I}_{g_{\gamma}},\quad \alpha \cdot
\mathcal{I}_{f_{\gamma}}=\mathcal{I}_{f_{\gamma}}\,\,\mbox{for}\,\,\alpha>0\quad\mbox{and}\quad
0\cdot \mathcal{I}_{f_{\gamma}}= \emptyset.
\end{equation*}
Suppose
$(\mathcal{Q},\mathcal{W})=(\sum_{\gamma\in\Gamma}\mathcal{P}_{\gamma},\times_{\gamma\in\Gamma}\varepsilon_{\mathcal{Y}_{\gamma}})$
and for each $\gamma\in\Gamma,$ consider the linear mapping
$\mathrm{T}_{\gamma}:\mathcal{P}_{\gamma}\to\texttt{Q}_{\gamma}$ for
all  $f_{\gamma}\in\mathcal{P}_{\gamma}$ by
$\mathrm{T}_{\gamma}(f_{\gamma})=\mathcal{I}_{f_{\gamma}}.$ Then
$\sum_{\gamma\in\Gamma}\mathrm{T}_{\gamma}(\mathcal{P}_{\gamma})=\sum_{\gamma\in\Gamma}\texttt{Q}_{\gamma}=\texttt{Q},$
where $\texttt{Q}=\{\sum_{\gamma\in\Delta}\mathcal{I}_{f_{\gamma}}:
f\in\mathcal{Q},f=\sum_{\gamma\in\Delta}f_{\gamma}\}.$ For elements
$\texttt{f},\texttt{g}\in\texttt{Q},$
$\texttt{f}=\sum_{\gamma\in\Delta}\mathrm{T}_{\gamma}(f_{\gamma}),$
$\texttt{g}=\sum_{\gamma\in\Theta}\mathrm{T}_{\gamma}(g_{\gamma}),$
the inductive limit order on $\texttt{Q}$ is defined by $
\texttt{f}\leq_{\overrightarrow{\Gamma}} \texttt{g}$ if $
\mathcal{I}_{f_{\gamma}}\subseteq
 \mathcal{I}_{g_{\gamma}}$ for all $\gamma\in\Delta\cup\Theta$
 and the inductive limit neighborhood for $w\in\mathcal{W},$
 $w=\times_{\gamma\in\Gamma}\varepsilon_{Y_{\gamma}}$ is
 defined by
 $\texttt{f}\leq_{\overrightarrow{\Gamma}}\texttt{g}+\overrightarrow{\times}_{\gamma\in\Gamma}\varepsilon_{Y_{\gamma}}$
 if $\mathcal{I}_{f_{\gamma}}\cap Y_{\gamma}\subseteq
 \mathcal{I}_{g_{\gamma}}$ for all $\gamma\in\Delta\cup\Theta.$
 Now, if
 we set $\mathcal{I}_{f}=\{x\in
X:f_{\gamma}(x_{\gamma})=+\infty\,\,\mbox{for
all}\,\,\gamma\in\Delta\}$ for all  $f\in\mathcal{Q},$
$f=\sum_{\gamma\in\Delta}f_{\gamma},$
$X=\times_{\gamma\in\Gamma}X_{\gamma}$ and
$\mathcal{Y}=\times_{\gamma\in\Gamma}\mathcal{Y}_{\gamma},$ then
$\texttt{Q}=\{\mathcal{I}_{f}:f\in\mathcal{Q}\}$ and $\mathcal{Y}$
is a family  subsets of $X$ which are downward directed with
respect to set inclusion. For
$Y\in\mathcal{Y},Y=\times_{\gamma\in\Gamma}Y_{\gamma}$ and
$\varepsilon>0,$ we define the neighborhood $\varepsilon_{Y}$ on
$\texttt{Q}$ by $\texttt{f}\leq \texttt{g}+\varepsilon_{Y},$ if
\begin{equation*}
\mathcal{I}_{f_{(\Delta\cup\Theta)}}\cap\times_{\gamma\in\Gamma}Y_{\gamma}
\subseteq \mathcal{I}_{g_{(\Delta\cup\Theta)}},
 \end{equation*}  where
$f_{(\Delta\cup\Theta)}=\sum_{\lambda\in\Delta\cup\Theta}f_{\lambda},$
$g_{(\Delta\cup\Theta)}=\sum_{\lambda\in\Delta\cup\Theta}g_{\lambda},$
$f_{\lambda}=f_{\gamma}$ for $\lambda\in\Delta,$
$g_{\lambda}=g_{\gamma}$ for $\lambda\in\Theta$ and
$f_{\lambda}=g_{\lambda}=0$ for $\lambda\in(\Delta\setminus
\Theta)\cup(\Theta \setminus\Delta).$ Then
$\varepsilon_{\mathcal{Y}}$ is identical to
$\overrightarrow{\mathcal{W}},$ where
$\varepsilon_{\mathcal{Y}}=\{\varepsilon_{Y}:Y\in\mathcal{Y}\}$ and
$\overrightarrow{\mathcal{W}}=\{\overrightarrow{w}:w\in\mathcal{W} ,
w=\times_{\gamma\in\Gamma}{\varepsilon_{Y}}_{\gamma}\},$
   that is, $(\texttt{Q},\varepsilon_{\mathcal{Y}})$ carries the locally convex
 inductive limit cone of
 $(\mathcal{P}_{\gamma},\varepsilon_{\mathcal{Y_{\gamma}}})$ by
  mappings $\mathrm{T}_{\gamma}.$
 \end{example}
\par  A \textit{dual pair} $(\mathcal{P},\mathcal{Q})$ consists of two
cones $\mathcal{P}$ and $\mathcal{Q}$ with a bilinear mapping
$(a,x)\longmapsto \langle a,x \rangle \,:\,\mathcal{P}\times
\mathcal{Q} \longrightarrow \overline{\mathbb{R}}.$ Let
$(\mathcal{P},\mathcal{Q})$ be a dual pair and $X$ be a collection
of subsets of $\mathcal{Q}$ such that:
\begin{enumerate}
    \item[(p$_{0}$)] $\inf \{\langle a,x \rangle : x \in A \} >-
    \infty$ for all $a\in \mathcal{P}$ and  $A \in X$.
    \item[(p$_{1}$)] $\lambda A \in X$ for all $A\in X$ and $\lambda >
    0$.
    \item[(p$_{2}$)] For all  $A,B\in X$ there is some $C\in
    X$
    such that $ A \cup B\subseteq C$.
\end{enumerate}
For each $A \in X$ we define $ \mathrm{U}_{A}=\{(a,b)\in
\mathcal{P}\times \mathcal{P}: \langle a,x\rangle \leq \langle b,x
\rangle +1 \,\, \mbox {for  all} \,\, x \in A\} .$ The set of all
$\mathrm{U}_{A}, A \in X$ forms a convex quasi-uniform structure
with property (U5) in \cite [Ch I, 5.2] {KWR} and defines a
locally convex structure on $\mathcal{P}$. This is called the
$X$\textit{-topology} on $\mathcal{P}$. For each $A\in X$ we
denote by $v_{A}$ the abstract neighborhood induced on
$\mathcal{P}$ by $\mathrm{U}_{A}$. Therefore $(a,b)\in
\mathrm{U}_{A}$ if and only if $a\leq b+v_{A}.$
\begin{remark} \label{Rem:2.13}
 (i) Every $X$-topology has (SP), also if
$X=\{v^{\circ}:v\in\mathcal{V}\},$ then $\mathcal{P}$ carries the
$X$-topology with respect to the dual pair
$(\mathcal{P},\mathcal{P}^{\ast});$ in particular, if
$(\mathcal{P},\mathcal{V})$ has the strict separation property then,
for each $v\in\mathcal{V},$ we have $\widetilde{v}\subseteq
\mathrm{U}_{v^{\circ}}\subseteq 2\/\widetilde{v}$ which implies that
$(\mathcal{P},\mathcal{V})$ is
 equivalent to $(\mathcal{P},\mathcal{V}_{X}).$ Note
that for an $X$-topology $(\mathcal{P},\mathcal{V}_{X}),$ we have
$\mathrm{U}_{A}=\mathrm{U}_{v^{\circ}_{A}}$ for all $A\in X$ so the
two abstract neighborhoods $\mathcal{V}_{X}$ and $\mathcal{V}_{X'}$
are identical to each other,
 where $X'=\{v^{\circ}_{A}:A\in X\}$ (cf. \cite[Remark 2
 (i)]{MRlcpl}).
 \par (ii)  If $X$ is a collection of subsets of $\mathcal{Q}$ satisfying
  $(p_{0}),(p_{1}),(p_{2})$ with respect to the dual pair
$(\mathcal{P},\mathcal{Q})$ then,  for every $A\in X,$ we have
$A\subseteq v^{\circ}_{A} = v^{\circ}_{v^{\circ}_{A}},$ where $
v^{\circ}_{v^{\circ}_{A}}=(v_{v^{\circ}_{A}})^{\circ};$ in
particular, if $(\mathcal{P},\mathcal{V})$ is a locally convex cone,
for each $v\in \mathcal{V},$ we have
$v^{\circ}=v^{\circ}_{v^{\circ}}.$ See \cite[Proposition 2.5]{MRHd}.
\end{remark}
\par The \emph{adjoint operator} of  
a linear operator $t:(\mathcal{P},\mathcal{V})\to
(\mathcal{R},\mathcal{Z})$ is defined by
\begin{eqnarray*}
\mu\to t'\mu: \mathcal{R}^{\ast}\to
\mathcal{L}(\mathcal{P},\overline{\mathbb{R}}),\quad t'\mu
(x)=(\mu\circ t)(x):\mathcal{P}\to\overline{\mathbb{R}},
\end{eqnarray*}
where $\mathcal{L}(\mathcal{P},\overline{\mathbb{R}})$ denotes the
linear mappings on $\mathcal{P}$ with values in
$\overline{\mathbb{R}}.$ According to Corollary 2.5 (a), we have
$\mathrm{T}'({\overrightarrow{w}}^{\circ})\subseteq w^{\circ}$ for
all $w\in\mathcal{W};$ in particular,
$\mathrm{T}'(\texttt{Q}^{\ast})\subseteq \mathcal{Q}^{\ast},$ where
$\mathrm{T}'$ is the adjoint operator of the mapping
$\mathrm{T}:(\mathcal{Q},\mathcal{W})\to(\texttt{Q},\overrightarrow{\mathcal{W}})$
(cf. \cite [Theorem 2.3 (a)]{MRHd}).
\begin{theorem}\label{Thm:2.14}
For each $\gamma\in\Gamma,$ let $X_{\gamma}$  be a collection of
subsets of $\mathcal{Q}_{\gamma}$ satisfying
$(p_{0}),(p_{1}),(p_{2})$ in dual pair
$(\mathcal{P}_{\gamma},\mathcal{Q}_{\gamma}).$ If
$(\texttt{Q},\overrightarrow{\mathcal{W}})$ is the locally convex
inductive limit cone of
$(\mathcal{P}_{\gamma},\mathcal{V}_{X_{\gamma}})$ by linear mappings
$\mathrm{T}_{\gamma},$
$(\mathcal{Q},\mathcal{W})=(\sum_{\gamma\in\Gamma}\mathcal{P}_{\gamma},\times_{\gamma\in\Gamma}\mathcal{V}_{X_{\gamma}}),$
 then
\begin{enumerate}
\item [$(\mathrm{a})$]  $\texttt{Q}$ has the $X$-topology with respect to the dual pair $(\texttt{Q},\texttt{Q}^{\ast}),$ where
$X=\{{\overrightarrow{w}}^{\circ}:w\in\mathcal{W},w=\times_{\gamma\in\Gamma}v_{A_{\gamma}}\},$
\item [$(\mathrm{b})$]  $\mathcal{Q}$ has the
$
\mathrm{T}'(X)$-topology with respect to the dual pair
$(\mathcal{Q},\mathcal{Q}^{\ast}),$
\item [$(\mathrm{c})$]
$(\mathcal{Q},\mathcal{V}_{\mathrm{T}'(X)})$  may be embedded on
$(\texttt{Q},\overrightarrow{\mathcal{W}}).$
 \end{enumerate}
\end{theorem}
\begin{proof}
Part (a)  holds by Remark \ref{Rem:2.13} (i). For (b), if
$a\in\mathcal{Q}$ and $w\in\mathcal{W},$ then
$0\leq_{\overrightarrow{\Gamma}} \mathrm{T}(a)+\lambda
\overrightarrow{w}$ for some $\lambda>0$ which yields $\inf\{\langle
a,\mathrm{T}'(\mu)\rangle:\mu\in
{\overrightarrow{w}}^{\circ}\}=\inf\{\langle
\mathrm{T}(a),\mu\rangle:\mu\in
{\overrightarrow{w}}^{\circ}\}\geq-\lambda>-\infty,$ i.e., $(p_{0})$
holds.  The properties $(p_{1}),(p_{2})$ are clear. For (c), we have
$a\leq b+v_{\mathrm{T}'({\overrightarrow{w}}^{\circ})}$ for
$a,b\in\mathcal{Q}$ and $w\in\mathcal{W}$ if and only if
$$\mu(\texttt{a})=(\mu\circ\mathrm{T})(a)=\langle a,\mathrm{T}'\mu\rangle\leq\langle
b,\mathrm{T}'\mu\rangle+1= (\mu\circ\mathrm{T})(b)
+1=\mu(\texttt{b})+1$$ for all $\mu\in
{\overrightarrow{w}}^{\circ},$ i.e., $\texttt{a}\leq
\texttt{b}+v_{{\overrightarrow{w}}^{\circ}}.$ That is,
$a\to\texttt{a}:(\mathcal{Q},\mathcal{V}_{\mathrm{T}'(X)})\to(\texttt{Q},\mathcal{V}_{X})$
is an embedding.
\end{proof}
\begin{corollary} \label{Cor:2.15}
 Let
 $(\texttt{Q},\overrightarrow{\mathcal{W}})$ be the
locally convex inductive limit cone of
$(\mathcal{P}_{\gamma},\mathcal{V}_{\gamma})$ by the linear
mappings $\mathrm{T}_{\gamma}$ and
$(\mathcal{Q},\mathcal{W})=(\sum_{\gamma\in\Gamma}\mathcal{P}_{\gamma},\times_{\gamma\in\Gamma}\mathcal{V}_{\gamma}).$
If $X_{\gamma}=\{v^{\circ}_{\gamma}: v_{\gamma}\in
\mathcal{V}_{\gamma}\}$ for all $\gamma\in\Gamma,$
$(\mathcal{Q},\mathcal{W}')=(\sum_{\gamma\in\Gamma}\mathcal{P}_{\gamma},\times_{\gamma\in\Gamma}\mathcal{V}_{X_{\gamma}})$
and  $(\texttt{Q},\overrightarrow{\mathcal{W}'})$ is the locally
convex inductive limit cone of
$(\mathcal{P}_{\gamma},\mathcal{V}_{X_{\gamma}})$ by
$\mathrm{T}_{\gamma},$
 then
\begin{enumerate}
\item [$(\mathrm{a})$]  $\texttt{Q}$ has the $X,Y$-topologies with respect to the dual pairs $(\texttt{Q},\texttt{Q}_{\mathcal{W}}^{\ast})$ and
$(\texttt{Q},\texttt{Q}_{\mathcal{W}'}^{\ast});$ respectively,
where
$X=\{{\overrightarrow{w}}^{\circ}:w\in\mathcal{W},w=\times_{\gamma\in\Gamma}v_{\gamma}\}$
and
$Y=\{{\overrightarrow{w'}}^{\circ}:w'\in\mathcal{W}',w'=\times_{\gamma\in\Gamma}v_{v^{\circ}_{\gamma}}\},$
\item [$(\mathrm{b})$]  $\mathcal{Q}$ has the
$\mathrm{T}'(X),\mathrm{T}'(Y)$-topologies with respect to the dual
pairs $(\mathcal{Q},\mathcal{Q}_{\mathcal{W}}^{\ast})$ and
$(\mathcal{Q},\mathcal{Q}_{\mathcal{W'}}^{\ast});$ respectively,
\item[$(\mathrm{c})$]
$(\mathcal{Q},\mathcal{V}_{\mathrm{T}'(Y)}$ may be embedded on
$(\texttt{Q},\mathcal{V}_{Y}),$
\item [$(\mathrm{d})$] if the all $(\mathcal{P}_{\gamma},\mathcal{V}_{\gamma})$
have (SP), then $(\mathcal{Q},\mathcal{V}_{\mathrm{T}'(X)})$ may be
embedded on $(\texttt{Q},\overrightarrow{\mathcal{W}}).$
 \end{enumerate}
\end{corollary}
\begin{proof}
Parts (a), (b) and (c) hold  by Theorem \ref{Thm:2.14}. For (d), by
Remark \ref{Rem:2.12} and Remark \ref{Rem:2.13} (ii), we have
  $w^{\circ}
  ={w'}^{\circ}\quad\mbox{for all}\quad  w\in
  \mathcal{W},w=\times_{\gamma\in\Gamma}v_{\gamma}.$
  Since $\mathrm{T'}\mu=\mu\circ \mathrm{T},$  by Corollary \ref{Cor:2.10} (a), we have $\mathrm{T}'\mu\in w^{\circ}$ if and only if
$\mathrm{T}'\mu\in {w'}^{\circ},$ that is,
${\overrightarrow{w}}^{\circ}={\overrightarrow{w'}}^{\circ}$ and
$\mathrm{T}'({\overrightarrow{w}}^{\circ})=\mathrm{T}'({\overrightarrow{w'}}^{\circ}).$
Thus, $(\texttt{Q},\mathcal{V}_{Y})$ is equivalent to
$(\texttt{Q},\mathcal{V}_{X})$ and
$(\mathcal{Q},\mathcal{V}_{\mathrm{T}'(Y)})$ is equivalent to
$(\mathcal{Q},\mathcal{V}_{\mathrm{T}'(X)}).$ On the other hand, by
Corollary \ref{Cor:2.11},
$(\texttt{Q},\overrightarrow{\mathcal{W}})$ has (SP), so
$(\texttt{Q},\mathcal{V}_{Y})$ is equivalent to
$(\texttt{Q},\overrightarrow{\mathcal{W}})$  by Remark
\ref{Rem:2.13} (i). Consequently,
$(\mathcal{Q},\mathcal{V}_{\mathrm{T}'(X)})$ may be embedded on
$(\texttt{Q},\overrightarrow{\mathcal{W}})$  by Theorem
\ref{Thm:2.14} (c).
\end{proof}
\begin{theorem}\label{Thm:2.16}
Let $(\texttt{Q},\overrightarrow{\mathcal{W}})$ be the locally
convex inductive limit cone of
$(\mathcal{P}_{\gamma},\mathcal{V}_{\gamma})$ by the linear mappings
$\mathrm{T}_{\gamma},$
$(\mathcal{Q},\mathcal{W})=(\sum_{\gamma\in\Gamma}\mathcal{P}_{\gamma},\times_{\gamma\in\Gamma}\mathcal{V}_{\gamma})$
and $X$ be a collection of subsets of $\mathcal{Q}$ satisfying
$(p_{0}),(p_{1}),(p_{2})$ with respect to the dual pair
$(\mathcal{Q}^{\ast},\mathcal{Q}).$ Then
\begin{enumerate}
\item [$(\mathrm{a})$] $\mathrm{T}(X)$ has
$(p_{0}),(p_{1}),(p_{2})$ with respect to the  dual pair
$(\texttt{Q}^{\ast},\texttt{Q}),$
    \item[$(\mathrm{b})$]
  $(\texttt{Q}^{\ast},\mathcal{V}_{\mathrm{T}(X)})$ may be embedded
in  $(\mathcal{Q}^{\ast},\mathcal{V}_{X}).$
\end{enumerate}
\end{theorem}
\begin{proof} (a) If $A\in X$  and $\mu\in\texttt{Q}^{\ast},$ then $
\inf\{\langle \mathrm{T}(a),\mu\rangle:a\in A\}=\inf\{\langle
a,\mathrm{T}' \mu\rangle:a\in A\}>-\infty$ so $(p_{0})$ holds. The
properties $(p_{1}),$ $(p_{2})$ are clear. For (b), if
$\mu,\nu\in\texttt{Q}^{\ast}$ and $A\in X,$ then $\mu\leq
\nu+v_{\mathrm{T}(A)}$ if and only if 
$$\langle \mathrm{T}'\mu,a
\rangle=(\mu\circ\mathrm{T})(a)\leq(\mu\circ\mathrm{T})(b)+1=\langle
\mathrm{T}'\nu,a \rangle+1$$ for all $a\in A,$ i.e.,
$\mathrm{T}'\mu\leq \mathrm{T}'\nu+v_{A}.$ That is,
$\mu\to\mathrm{T}'\mu:(\texttt{Q}^{\ast},\mathcal{V}_{\mathrm{T}(X)})\to(\mathcal{Q}^{\ast},\mathcal{V}_{X})$
is an embedding.
\end{proof}

\par  If $(\mathcal{P},\mathcal{V})$ is a locally convex cone then, the resulting $X$-topology
 on $\mathcal{P}^{\ast}$ with respect to the dual pair $(\mathcal{P}^{\ast},\mathcal{P})$
 is called the \emph{weak*-topology} and denoted by
 $\sigma(\mathcal{P}^{\ast},\mathcal{P}).$ The \emph{weak topology}
 $\sigma(\mathcal{P},\mathcal{P}^{\ast})$ is defined in the similar way.
 Also,   $w(\mathcal{P}^{\ast},\mathcal{P})$  denotes the \emph{topology
of pointwise convergence }of the elements of $\mathcal{P}$
considered as functions on $\mathcal{P}^{\ast}$ with values in
$\overline{\mathbb{R}}.$ A typical neighborhood for an element
$\nu\in\mathcal{P}^{\ast}$ in $w(\mathcal{P}^{\ast},\mathcal{P})$
is given via a finite subset
 $A=\{a_{1},a_{2},...,a_{n}\}$ of $\mathcal{P}$ by
\begin{eqnarray*}
w_{A}(\nu)&=&\left\{%
\mu\in\mathcal{P}^{\ast}:
\begin{array}{ll}
|\mu(a_{i})-\nu(a_{i})|\leq1,
\vspace{1mm} &\hbox{if $\nu(a_{i})<+\infty$}\\
\mu (a_{i})\geq 1\vspace{0mm}, &\hbox{if $\nu(a_{i})=+\infty$}
\end{array}%
\right\}.
\end{eqnarray*}
For details see [1, Ch II, 3.4]. If we denote by
$\tau(\mathcal{P},\mathcal{P}^{\ast})$ the topology of all
$w(\mathcal{P}^{\ast},\mathcal{P})$-compact subsets of
$\mathcal{P}^{\ast}$ on $\mathcal{P}$ then, for each
$v\in\mathcal{V},$  $v^{\circ}$ is
$w(\mathcal{P}^{\ast},\mathcal{P})$-compact by [1, Ch II,
Proposition 2.4], so
$(\mathcal{P},\mathcal{V}_{\tau(\mathcal{P},\mathcal{P}^{\ast})})$
is finer than the $X$-topology $(\mathcal{P},\mathcal{V}_{X}),$
where $X=\{v^{\circ}:v\in \mathcal{V}\}.$ Thus,
$(\mathcal{P},\mathcal{V}_{\tau(\mathcal{P},\mathcal{P}^{\ast})})$
will be finer than the original topology
$(\mathcal{P},\mathcal{V}),$ whenever $(\mathcal{P},\mathcal{V})$
has the strict separation property (see \cite [Proposition
3.7]{MRlcpd}).
\begin{corollary}
If $(\texttt{Q},\overrightarrow{\mathcal{W}})$ is the locally convex
inductive limit  cone of
  $(\mathcal{P}_{\gamma},\mathcal{V}_{\gamma})$ by the linear mappings
  $\mathrm{T}_{\gamma}$ and $(\mathcal{Q},\mathcal{W})=
     (\sum_{\gamma\in\Gamma}\mathcal{P}_{\gamma},\times_{\gamma\in\Gamma}\mathcal{V}_{\gamma}),$
          then
\begin{enumerate}
\item [$(\mathrm{a})$]  $(\texttt{Q}^{\ast},\mathcal{V}_{\sigma(\texttt{Q}^{\ast},\texttt{Q})})$ may be embedded
in
$(\mathcal{Q}^{\ast},\mathcal{V}_{\sigma(\mathcal{Q}^{\ast},\mathcal{Q})}),$
\item[$(\mathrm{b})$] $\mathrm{T}'$ is continuous with respect to
$w(\texttt{Q}^{\ast},\texttt{Q})$ on $\texttt{Q}^{\ast}$ and
$w(\mathcal{Q}^{\ast},\mathcal{Q})$ on $\mathcal{Q}^{\ast},$
\item [$(\mathrm{c})$]
$(\mathcal{Q},\mathcal{V}_{\tau(\mathcal{Q},\mathcal{Q}^{\ast})})$
may be embedded on
$(\texttt{Q},\mathcal{V}_{\tau(\texttt{Q},\texttt{Q}^{\ast})}),$
\item [$(\mathrm{d})$]
$(\texttt{Q},\mathcal{V}_{\tau(\texttt{Q},\texttt{Q}^{\ast})})$ is
the locally convex inductive limit cone of
$(\mathcal{P}_{\gamma},\mathcal{V}_{\tau(\mathcal{P}_{\gamma},\mathcal{P}_{\gamma}^{\ast})})$
by mappings $\mathrm{T}_{\gamma}.$
\end{enumerate}
\end{corollary}
\begin{proof} Part (a) holds by Theorem \ref {Thm:2.16} (b). For (b), we show
that for each $\nu\in\mathcal{Q}^{\ast}$ and
$A\subset\mathcal{Q},\,A=\{a_{1},a_{2},...,a_{n}\},$
$\mathrm{T}'(w_{\mathrm{T}(A)}(\nu))\subseteq
w_{A}(\mathrm{T}'(\nu)).$ Let $\mu\in w_{\mathrm{T}(A)}(\nu).$ If
$\mathrm{T}'\nu (a_{i})=(\nu\circ\mathrm{T})(a_{i})<+\infty,$ then
$|\mathrm{T}'\nu(a_{i})-\mathrm{T}'\mu(a_{i})|=|(\nu\circ\mathrm{T})(a_{i})-(\mu\circ\mathrm{T})(a_{i})|\leq1,$
and if $\mathrm{T}'\nu
(a_{i})=(\nu\circ\mathrm{T})(a_{i})=+\infty,$ then
$\mathrm{T}'\mu(a_{i})=(\mu\circ \mathrm{T})(a_{i})\geq 1.$ That
is $\mathrm{T}'\mu\in w_{A}(\mathrm{T}'(\nu)).$
 For (c),  if $\Upsilon\subset
\texttt{Q}^{\ast}$ is $w(\texttt{Q}^{\ast},\texttt{Q})$-compact,
then $\mathrm{T}'(\Upsilon)\subseteq \mathcal{Q}^{\ast}$ will be
$w(\mathcal{Q}^{\ast},\mathcal{Q})$-compact by (b). Thus, by Theorem
\ref{Thm:2.14} (c), we have $a\leq b+v_{\mathrm{T}'(\Upsilon)}$ for
$a,b\in\mathcal{Q}$ if and only if $\mathrm{T}(a)\leq
\mathrm{T}(b)+v_{\Upsilon},$ that is,
$a\to\mathrm{T}(a):(\mathcal{Q},\mathcal{V}_{\tau(\mathcal{Q},\mathcal{Q}^{\ast})})\to(\texttt{Q},\mathcal{V}_{\tau(\texttt{Q},\texttt{Q}^{\ast})})$
is an embedding. For (d), let
$(\texttt{Q},\overrightarrow{\mathcal{W}})$ be the locally convex
inductive limit cone of
$(\mathcal{P}_{\gamma},\mathcal{V}_{\tau(\mathcal{P}_{\gamma},\mathcal{P}^{\ast}_{\gamma})})$
by $\mathrm{T}_{\gamma}.$ According to \cite [Corollary 3.9]
{MRlcpd},
$(\mathcal{Q},\mathcal{V}_{\tau(\mathcal{Q},\mathcal{Q}^{\ast}})=(\sum_{\gamma\in\Gamma}
\mathcal{P}_{\gamma},\times_{\gamma\in\Gamma}\mathcal{V}_{\tau(\mathcal{P}_{\gamma},\mathcal{P}^{\ast}_{\gamma})}),$
so  by (b) and Proposition \ref{prop:2.7} (a), we infer that
$(\texttt{Q},\overrightarrow{\mathcal{W}})$ is finer than
$(\texttt{Q},\mathcal{V}_{\tau(\texttt{Q},\texttt{Q}^{\ast})}).$ On
the other hand, each polar ${\overrightarrow{w}}^{\circ}$ is
$w(\texttt{Q}^{\ast},\texttt{Q})$-compact by \cite [Ch II,
2.4]{KWR}, so
$(\texttt{Q},\mathcal{V}_{\tau(\texttt{Q},\texttt{Q}^{\ast})})$ is
identical to  $(\texttt{Q},\overrightarrow{\mathcal{W}}).$
\end{proof}

\bibliographystyle{plain}

\begin{thebibliography}{99}
\bibitem [1]{KWR}
Keimel, K., Roth W.: \textit{Ordered cones and approximation},
Lecture Notes in Mathematics, vol. 1517, Springer,
Heidelberg-Berlin-New York (1992)
\bibitem [2]{MRc} Motallebi, M.R.: \textit{Completeness on locally convex cones}, C. R.
Math. Acad. Sci.  Paris 352 (10),  785-789 (2014)
\bibitem [3]{MRlcpd} Motallebi, M.R.: \textit{Locally convex product and direct sum
cones}, Mediterr. J. Math. 11 (3),  913-927 (2014)
\bibitem [4]{MRlcpl} Motallebi, M.R.: \textit{Locally convex projective limit
cones}, Math. Slovaca 66 (6), 1387-1398 (2016)
\bibitem [5]{MRo} Motallebi, M.R.: \textit{On weak completeness of products and direct sums in locally convex
cones}. Period. Math. Hung. (2017). doi: 10.1007/s10998-017-0201-4
\bibitem [6]{MRHd} Motallebi, M.R.,
 Saiflu. H.: \textit{Duality on locally convex cones}, J. Math. Anal.
Appl. 337 (2),  888-905 (2008)
\bibitem [7]{MRHpd} Motallebi, M.R., Saiflu. H.: \textit{Products and direct sums in
locally convex cones}, Can. Math. Bull. 55 (4),  783-798 (2012)
\bibitem [8]{WRl}
Roth, W.: \textit{Locally convex quotient cones}, J. Convex Anal.
18 (4), 903-913 (2011)
\bibitem [9]{WRo} Roth, W.: \textit{Operator-valued measures and integrals for
cone-valued functions}, Lecture Notes in Mathematics, vol. 1964,
Springer Verlag, Heidelberg-Berlin-New York (2009)
\end{thebibliography}

\end{document}